# Histogram Arithmetic under Uncertainty of Probability Density Function

**Vladimir N. Petrushin**

Moscow State University of Printing Arts
Moscow, 127550 Russia

**Evgeny V. Nikulchev**

Moscow Technological Institute
Leninskiy pr., 38A, Moscow, 119334 Russia

**Dmitry A. Korolev**

Moscow State Technical University of Radio Engineering, Electronics and Automatics, Moscow, 119454 Russia

**Abstract**

In this article we propose a method of performing arithmetic operations on variables with unknown distribution. The approach to the evaluation results of arithmetic operations can select probability intervals of the algebraic equations and their systems solutions, of differential equations and their systems in case of histogram evaluation of the empirical density distributions of random parameters.

## 1 Introduction

The solution of practical problems in various areas of scientific and practical activities forces us to find out the value of point solutions as well as different estimates of these solutions in the form of sets. The uncertainty of the initial informa-tion in physics, chemistry, biology, medicine, engineering, sociology, economics and other fields of human activity generates various evaluation systems.

There exist three categories: interval, fuzzy and random. The first two categories during the large amounts of linear and nonlinear transformations give rise to many solutions; still their boundaries are not acceptable in practice.

The tasks with random parameters (e.g. management tasks, planning tasks, etc.) can be of great interest. Thereby it is vital to study the behavior of variables with random character during the arithmetic operations applied to them.

Many scholars studied the issue of interval and histogram arithmetic [1–3]. Recently, interval arithmetic is used in conjunction with fuzzy sets [4, 5]. B. Dobronets and O. Popova studied systematic histogram approach within the arithmetic operations evaluation for parametric tasks with known probability density parameters [6]. The issue of arithmetic operations with random variables, for which on the basis of experimental data you can build an empirical membership function, remains unresolved.

Let us consider the task of arithmetic operations for parameters with the unknown distribution. As a result of experiments, the researchers obtain a set (sample) of discrete parameters values regardless of the discrete or continuous random variable.

To realize the behavior of a studied random variable it is essential to assess its probability density function. Traditionally, on the basis of experimental data we should build a bar chart [1–3].

Such histograms will be further the research object of this article. The article's objective is to find out an acceptable numerical evaluation of the arithmetic operations results with random variables.

## 2 Histogram plotting based on the original sample

These rather strict rules of histogram plotting are formulated by the authors for the first time and based on the reliability of the group average within all the bar graph semi segments and their almost significant difference. The system of inequalities follows from the above:

$$\begin{cases} \overline{X}_j - \delta_j > x_j^{(\min)}; \\ \overline{X}_j + \delta_j < x_j^{(\max)}, \end{cases}$$

where $\overline{X}_j$ — the mathematical expectation of the group average of $j$ half-segment; $x_j^{(\min)}$, $x_j^{(\max)}$ — the boundaries of $j$ half-segment; $\delta_j = t(\gamma_j, n_j) \cdot S_j \left(1 + q(\gamma_j, n_j)\right) / \sqrt{n_j}$ — statistical error in defining the average according to the Student distribution, where $t(\gamma_j, n_j)$ — the value of a t-test for a given reliability and the scope of the group $n_j$; $S_j$ — unbiased assessment of the intragroup variance in $j$ half-segment; $q(\gamma_j, n_j)$ — amendment to the standard quadratic deviation.

The obtained system of inequalities allows us to create inequalities allocation groups selection algorithm in a sampling variation number on the basis of the task with the confidence probability within $\gamma_j$ half-segments .

The solution to this inequalities system is an algorithm of a group selection. Moreover, as a result of solving the system, we obtain the group volume as well as the boundaries of the half-segment. The given approach significantly limits voluntarism in the histogram construction, though it doesn't eliminate it completely.

If the condition is not met for the last segment, it must be connected to the previous one. But the knowledge of the extreme values of the sample included in a half-segment leaves the possibility of a particular boundaries variation between the right extreme value in the previous half-segment and the extreme left in the current.

Obviously, taking the hypothesis on the group average independence means that the histogram reliability in a whole $\gamma$ will be the reliability product of all the group averages $\gamma_j$. Thus we get the first component of the histogram quality assessment $\gamma = \prod_{j=1}^{k} \gamma_j$.

Since our main objective is to assess the approximation quality of the empirical distribution function obtained from the piecewise linear function histogram, the authors propose to use the probability value of the type I error according to Kolmogorov criterion as a measure [7].

Thus, the quality of the bar chart can be determined by the value $Q = \gamma\alpha$.

In such constructions the sample representativeness is very important. Quantitative criteria of this value don't exist at the moment. The issue solved by the researcher dictates the quantitative requirements of representativeness itself.

We can name a huge amount of knowledge areas the researchers are interested in, where the sample with the volume equal to one is representative. In this case, the sample representativeness is determined by the bar graph reliability, i.e. the reliability of all the group averages within the half-intervals and the corresponding variances.

## 3 Arithmetic operations with histogram given variables

Let us consider the possibility of arithmetic operations on variables, which distribution is defined by the empirical histogram (piecewise constant probability density function). The direct convolution operation use in this case is impossible [8], therefore we suggest the geometric method.

Assume that we have the histograms of two random variables $X$ and $Y$, based on the regulations given above, at the same time $x \in \left[\underline{a};\overline{a}\right]$, $y \in \left[\underline{b};\overline{b}\right]$, and a histogram of the $X$ quantity consists of $k$ half-segments, the value of $Y$ consists of $m$ half-slots.

While carrying out arithmetic operations with the values given by the constant piecewise density distribution, the operations are conducted within all the rectangles formed by combinations of half-segments variables $X$ and $Y$.

If there are a lot of such values, it is necessary to consider parallelotopes of possible combinations between random variables, but given the independence of these variables, folding procedure can be performed sequentially.

Let us find out the distribution function of the sum $(Z = X + Y)$, the difference $(Z = X - Y)$, the product $(Z = XY)$, the quotient $(Z = Y / X)$. Suppose for definiteness $\left(\overline{b}_r - \underline{b}_r\right) > \left(\overline{a}_r - \underline{a}_r\right)$, $\underline{a}_1 = \underline{a} > 0$, $\underline{b}_1 = \underline{b} > 0$, $j$ and $r$ — numbers of half-segments. The values of the function distribution are equal to the area cut off from the rectangle with sides $\left(\overline{a}_j - \underline{a}_j\right)$, $\left(\overline{b}_j - \underline{b}_j\right)$, normalized to the area of this rectangle:

- for the sum — the line $y = z - x$ (Fig. 1);
- for the difference — the line $y = x - z$ (Fig. 2);
- for the product — the hyperbola $y = z / x$ (Fig. 3);
- for the quotient — the line $y = zx$ (Fig. 4).

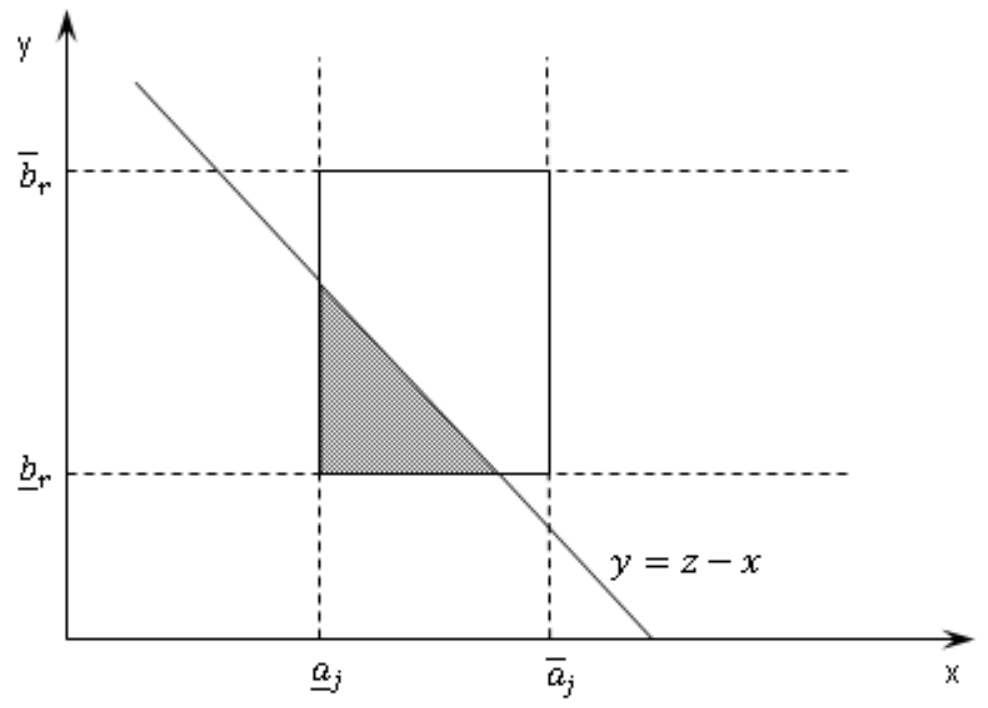

Fig. 1. The distribution function of the sum ( $Z = X + Y$ )

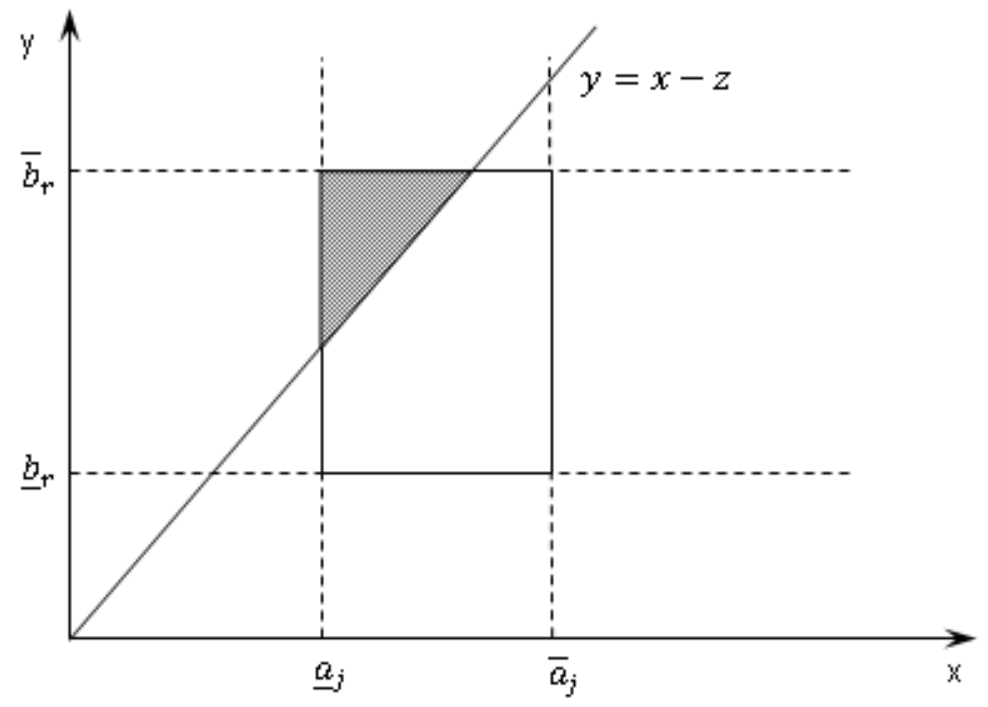

Fig. 2. The distribution function of the difference ( $Z = X - Y$ )

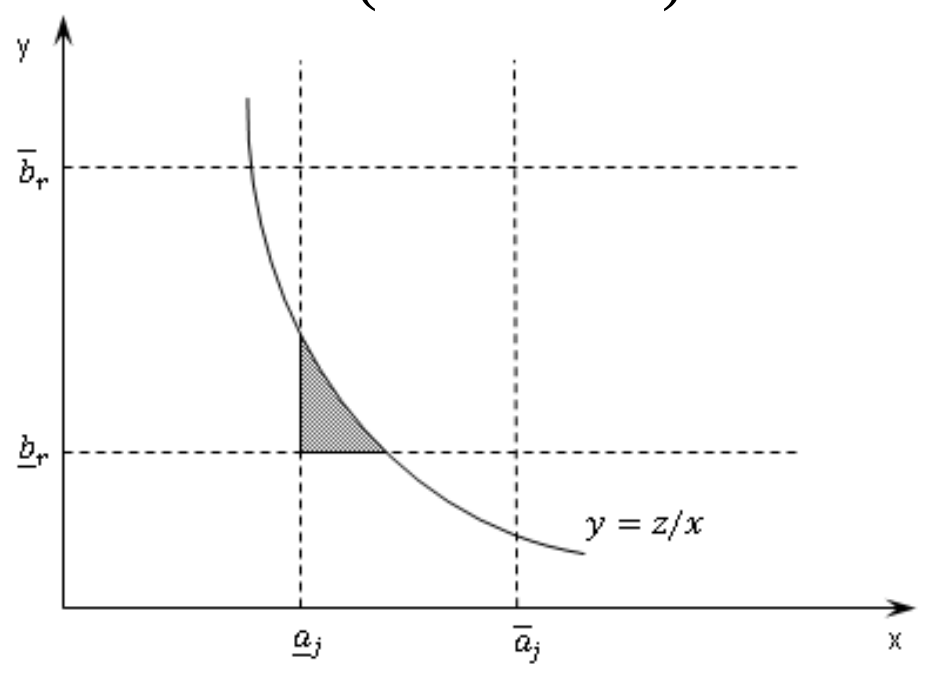

Fig. 3. The distribution function of the product ( $Z = XY$ )

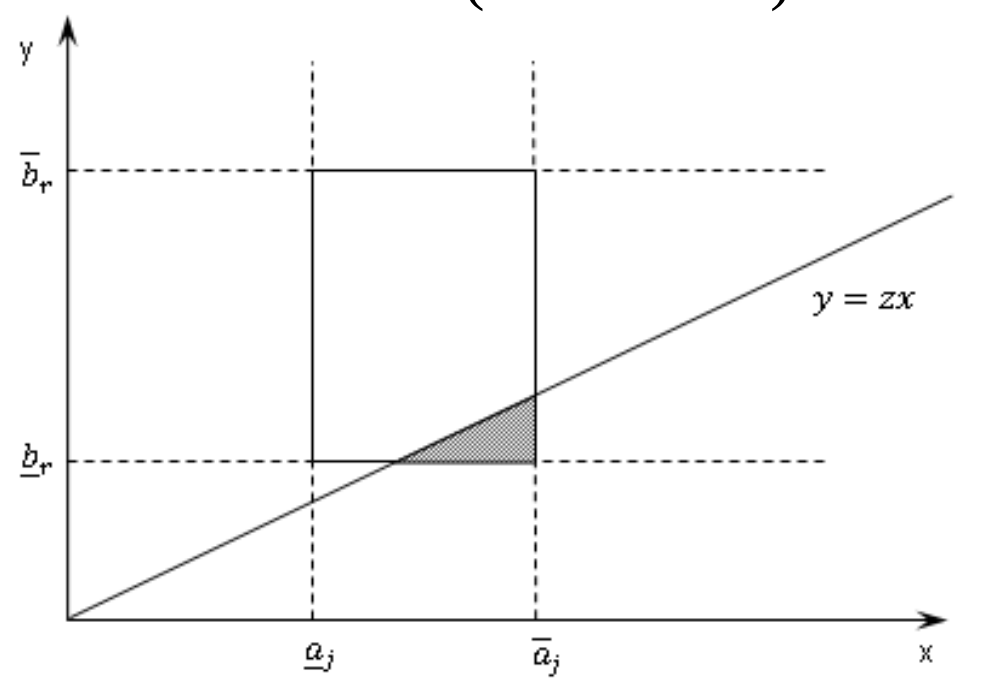

Fig. 4. The distribution function of the quotient ( $Z = Y / X$ )

The lines and the hyperbola are moving in the **grad** $Z$ direction. The distribution functions of random variables resulted from arithmetic operations.

Then the distribution function of the arithmetic operations results over the random variables within a specified rectangle is as follows:

$$F_{jr}\left(z=x+y\right)=\frac{1}{\left(\overline{a}-\underline{a}\right)\left(\overline{b}-\underline{b}\right)}\begin{cases}0, & z<\underline{a}_j+\underline{b}_r;\\ \frac{1}{2}\left(z-\underline{a}_j-\underline{b}_r\right)^2, & \underline{a}_j+\underline{b}_r\le z<\overline{a}_j+\underline{b}_r;\\ \frac{1}{2}\left(\overline{a}_j-\underline{a}_j\right)^2+\left(\overline{a}_j-\underline{a}_j\right)\left(z-\overline{a}_j-\underline{b}_r\right), & \overline{a}_j+\underline{b}_r\le z<\underline{a}_j+\overline{b}_r;\\ \left(\overline{a}_j-\underline{a}_j\right)\left(\overline{b}_r-\underline{b}_r\right)-\frac{1}{2}\left(\overline{a}_j-\overline{b}_r-z\right)^2, & \underline{a}_j+\overline{b}_r\le z\le\overline{a}_j+\overline{b}_r;\\ \left(\overline{a}_j-\underline{a}_j\right)\left(\overline{b}_r-\underline{b}_r\right), & z>\overline{a}_j+\overline{b}_r.\end{cases}$$

$$F_{jr}\left(z=x-y\right)=\frac{1}{\left(\overline{a}-\underline{a}\right)\left(\overline{b}-\underline{b}\right)}\begin{cases}0, & z<\underline{a}_j-\overline{b}_r;\\ 0.5\left(z-\underline{a}_j+\overline{b}_r\right)^2, & \underline{a}_j-\overline{b}_r\le z<\overline{a}_j-\overline{b}_r;\\ 0.5\left(\overline{a}_j-\underline{a}_j\right)^2+\left(\overline{a}_j-\underline{a}_j\right)\left(z-\overline{a}_j+\overline{b}_r\right), & \overline{a}_j-\overline{b}_r\le z<\underline{a}_j-\underline{b}_r;\\ \left(\overline{a}_j-\underline{a}_j\right)\left(\overline{b}_r-\underline{b}_r\right)-\frac{1}{2}\left(\overline{a}_j-\underline{b}_r-z\right)^2, & \underline{a}_j-\underline{b}_r\le z\le\overline{a}_j-\underline{b}_r;\\ \left(\overline{a}_j-\underline{a}_j\right)\left(\overline{b}_r-\underline{b}_r\right), & z>\overline{a}_j-\underline{b}_r.\end{cases}$$

$$F_{jr}\left(z=xy\right)=\frac{1}{\left(\overline{a}-\underline{a}\right)\left(\overline{b}-\underline{b}\right)}\begin{cases}0, & z<\underline{a}_j\underline{b}_r;\\ z\left(\ln z-\ln\underline{a}_j-\ln\underline{b}_r-1\right)^2, & \underline{a}_j\underline{b}_r\le z<\overline{a}_j\underline{b}_r;\\ z\left(\ln\overline{a}_j-\ln\underline{a}_j\right)^2-\underline{b}_r\left(\overline{a}_j-\underline{a}_j\right), & \overline{a}_j\underline{b}_r\le z<\underline{a}_j\overline{b}_r;\\ \left(\overline{a}_j-\underline{a}_j\right)\left(\overline{b}_r-\underline{b}_r\right)-\overline{a}_j\overline{b}_r-z\left(\ln z-\ln\overline{a}_j-\ln\underline{b}_r-1\right)^2, & \underline{a}_j\overline{b}_r\le z\le\overline{a}_j\overline{b}_r;\\ \left(\overline{a}_j-\underline{a}_j\right)\left(\overline{b}_r-\underline{b}_r\right), & z>\overline{a}_j\overline{b}_r.\end{cases}$$

$$F_{jr}\left(z=y/x\right)=\frac{1}{\left(\overline{a}-\underline{a}\right)\left(\overline{b}-\underline{b}\right)}\begin{cases}0, & z<\underline{b}_r/\overline{a}_j;\\ 0.5\left(\overline{a}_j\sqrt{z}-\frac{\underline{b}_r}{\sqrt{z}}\right)^2, & \underline{b}_r/\overline{a}_j\le z<\underline{b}_r/\underline{a}_j;\\ 0.5\left(\overline{a}_j-\underline{a}_j\right)^2 z, & \underline{b}_r/\underline{a}_j\le z<\overline{b}_r/\overline{a}_j;\\ \left(\overline{a}_j-\underline{a}_j\right)\left(\overline{b}_r-\underline{b}_r\right)-0.5\left(\underline{a}_j\sqrt{z}-\frac{\overline{b}_r}{\sqrt{z}}\right)^2, & \overline{b}_r/\overline{a}_j\le z\le\overline{b}_r/\underline{a}_j;\\ \left(\overline{a}_j-\underline{a}_j\right)\left(\overline{b}_r-\underline{b}_r\right), & z>\overline{b}_r/\underline{a}_j.\end{cases}$$

While changing the ratio of the sides lengths and the rectangle vertex coordinates, the form of the distribution function remains, the constants and intervals will change.

The distribution function within the full range of changing the $Z$ variable $F(z)=\sum_{j=1}^{k}\sum_{r=1}^{m}F_{jr}$ , and the probability density function

$$f(z)=F'(z)=\sum_{j=1}^{k}\sum_{r=1}^{m}F'_{jr}(z)=\sum_{j=1}^{k}\sum_{r=1}^{m}f_{jr}(z).$$

Differentiating the distribution functions of the arithmetic operations results over the random variables obtained above, we get:

$$f_{jr}(z=x+y)=\frac{1}{(\overline{a}-\underline{a})(\overline{b}-\underline{b})}\begin{cases}0, & z<\underline{a}_j+\underline{b}_r;\\ z-\underline{a}_j-\underline{b}_r, & \underline{a}_j+\underline{b}_r\le z<\overline{a}_j+\underline{b}_r;\\ \overline{a}_j-\underline{a}_j, & \overline{a}_j+\underline{b}_r\le z<\underline{a}_j+\overline{b}_r;\\ \overline{a}_j+\overline{b}_r-z, & \underline{a}_j+\overline{b}_r\le z\le\overline{a}_j+\overline{b}_r;\\ 0, & z>\overline{a}_j+\overline{b}_r.\end{cases}$$

$$f_{jr}(z=x-y)=\frac{1}{(\overline{a}-\underline{a})(\overline{b}-\underline{b})}\begin{cases}0, & z<\underline{a}_j-\overline{b}_r;\\ z-\underline{a}_j+\overline{b}_r, & \underline{a}_j-\overline{b}_r\le z<\overline{a}_j-\overline{b}_r;\\ \overline{a}_j-\underline{a}_j, & \overline{a}_j-\overline{b}_r\le z<\underline{a}_j-\underline{b}_r;\\ \overline{a}_j-\underline{b}_r-z, & \underline{a}_j-\underline{b}_r\le z\le\overline{a}_j-\underline{b}_r;\\ 0, & z>\overline{a}_j-\underline{b}_r.\end{cases}$$

$$f_{jr}(z=xy)=\frac{1}{(\overline{a}-\underline{a})(\overline{b}-\underline{b})}\begin{cases}0, & z<\underline{a}_j\underline{b}_r;\\ \ln^2 z-\ln^2(\underline{a}_j\underline{b}_r)-1, & \underline{a}_j\underline{b}_r\le z<\overline{a}_j\underline{b}_r;\\ (\ln\overline{a}_j-\ln\underline{a}_j)^2, & \overline{a}_j\underline{b}_r\le z<\underline{a}_j\overline{b}_r;\\ 1+\ln^2\overline{a}_j\overline{b}_r-\ln^2 z, & \underline{a}_j\overline{b}_r\le z\le\overline{a}_j\overline{b}_r;\\ 0, & z>\overline{a}_j\overline{b}_r.\end{cases}$$

$$f_{jr}(z=y/x)=\frac{1}{(\overline{a}-\underline{a})(\overline{b}-\underline{b})}\begin{cases}0, & z<\underline{b}_r/\overline{a}_j;\\ 0.5\left(\overline{a}_j^{2}-(\underline{b}_r^{2}/z^{2})\right), & \underline{b}_r/\overline{a}_j\le z<\underline{b}_r/\underline{a}_j;\\ 0.5(\overline{a}_j-\underline{a}_j)^2, & \underline{b}_r/\underline{a}_j\le z<\overline{b}_r/\overline{a}_j;\\ 0.5\left((\overline{b}_r^{\,2}/z^2)-\underline{a}_j^{\,2}\right), & \overline{b}_r/\overline{a}_j\le z\le\overline{b}_r/\underline{a}_j;\\ 0, & z>\overline{b}_r/\underline{a}_j.\end{cases}$$

Expressions for the probability distribution of the arithmetic operations results look not very bulky and they can be formally used while solving various tasks. However, if there are a lot of such actions, including non-linear combinations, the task becomes computationally complex and hard to solve analytically.

In addition, it's vital to note that the number of half-intervals, even when there is only one action, is equal to $4km$, where $k$ — the number of histogram $X$ segments, $m$ — the number of histogram $Y$ segments. Then there immediately appears the question whether these results in each of these intervals are reliable or not. Since the parameters values included in them may fail to provide reliable assessment of the group average.

Moreover, to solve such tasks by numerical methods we need specific software.

## 4 Statistical assessment of the histograms results

The result of any finite number of linear and non-linear transformations, including arithmetic operations, can be assessed statistically on the basis of discrete sampling.

This requires computing all the operations with all the possible combinations of the observed values of a random variable, and then plotting a histogram.

The half-interval reliability of such bar charts determines by the product of all the reliable assessments of group averages taken from random variables. Thus, the new interval will include the number of quantities $z$ that will ensure the reliability of each item, factor and so on.

Let us illustrate this with the specific example of two samples (samples adress). Relying on the rules above, let us build their histograms. For this let us set $\gamma_j = \gamma_r = 0,999$.

Based on the principle of group averages reliability with the picked reliability, we get the following histograms of random variables $X$ and $Y$ (Fig. 5, 6).

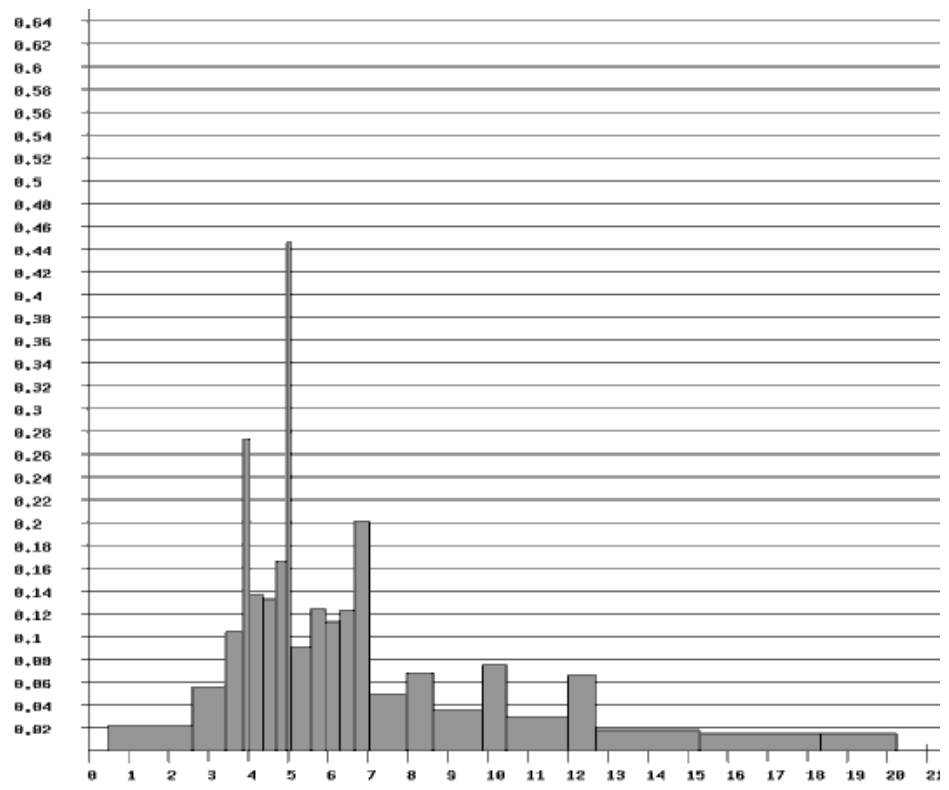

Fig. 5. Bar graph $X$

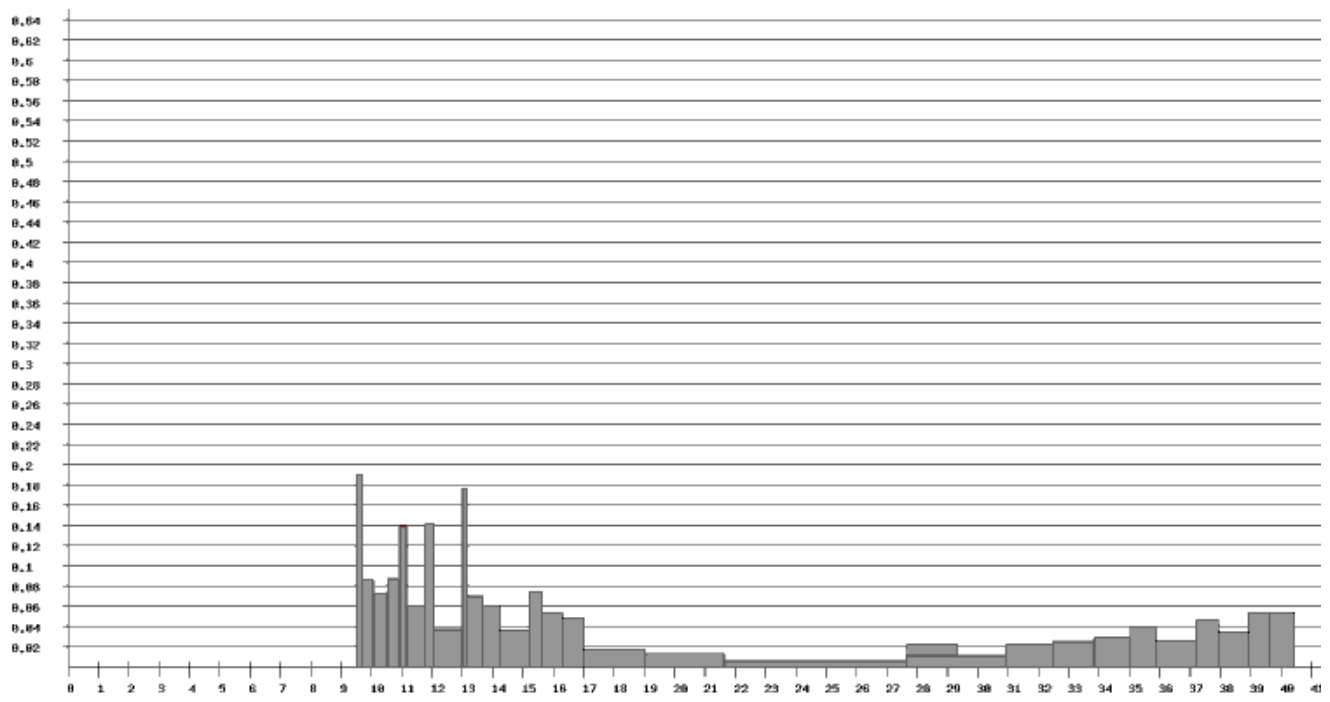

Fig. 6. Bar graph $Y$

Their construction was based on the averages assessments reliability for each half-segment of the arithmetic operation result. It means that the creation of the quantity $Z$ segment involves such amounts of random variables $X$ and $Y$ which guarantee selected reliability of the quantities accuracy.

We do not observe the quantity $Z$ as a result of actions over $X$ and $Y$. Thus, the half-segment reliability of $Z$ values is equal to the reliability product of the half-segments $X$ and $Y$, which formed the half-interval $Z$.

While constructing the histogram distribution of $X$ and $Y$ quantities, we checked all the group averages for accuracy within the histogram half-segments. If the question of arithmetic operations results approximation is vital for a researcher, it is appropriate to do it after rather than calculate the convolution of approximated summands and factors.

Fig. 7–10 represent bar graphs obtained after the arithmetic operations with samples $X$ and $Y$.

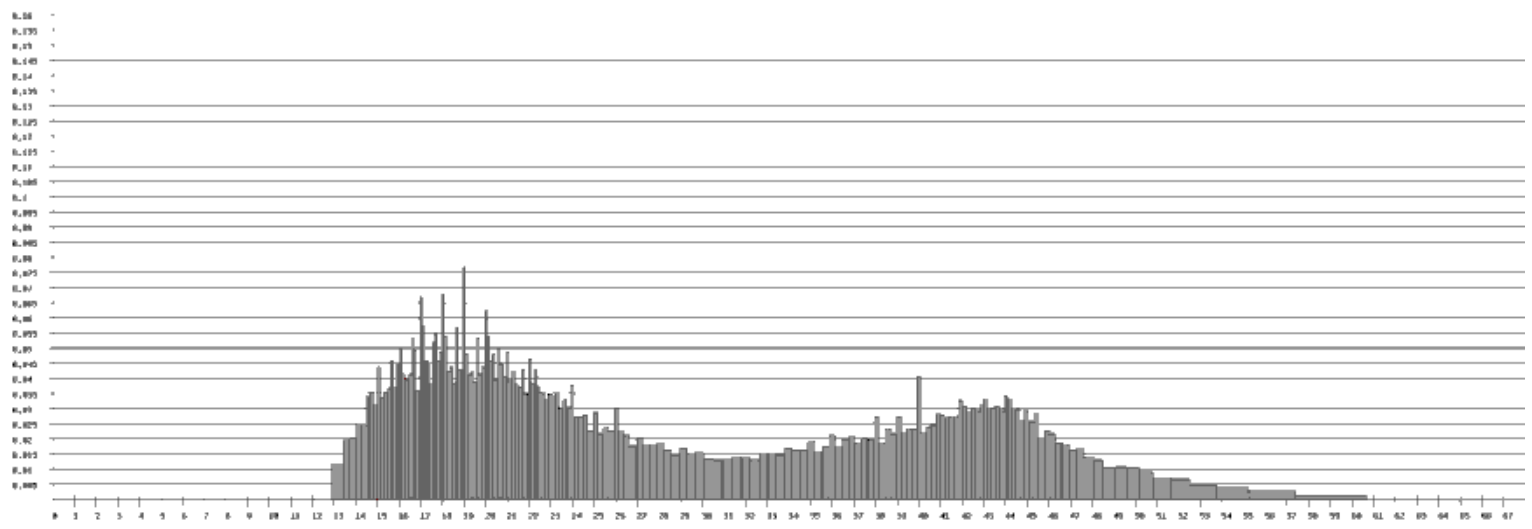

Fig. 7. Bar graph $Z = X + Y$

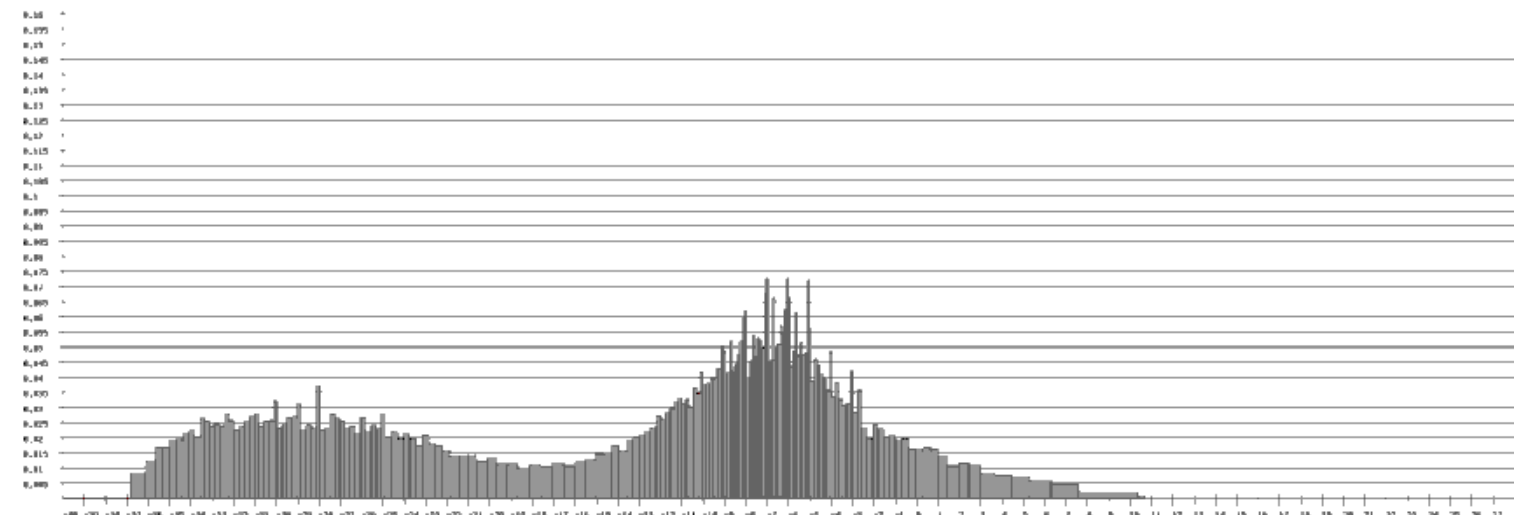

Fig. 8. Bar graph $Z = X - Y$

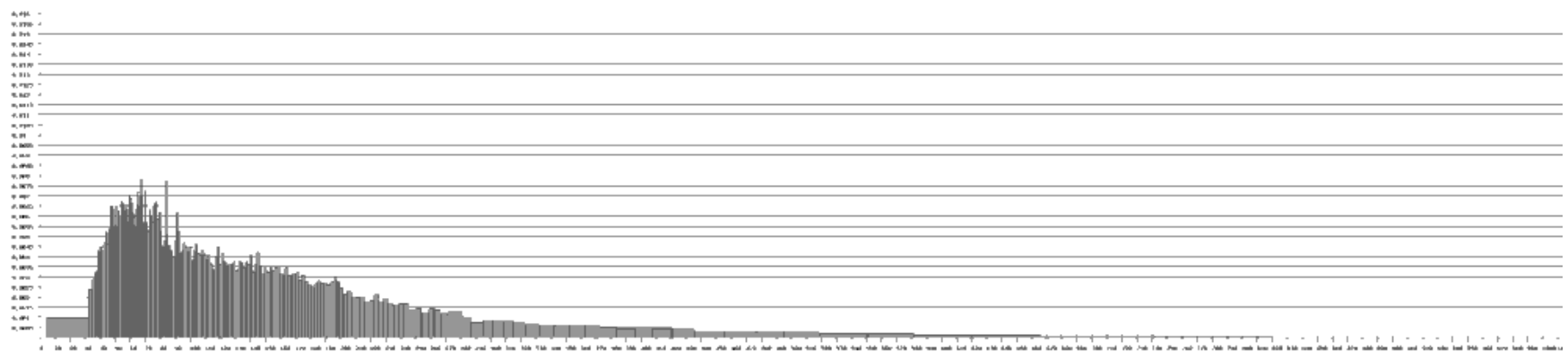

Fig. 9. Bar graph $Z = XY$

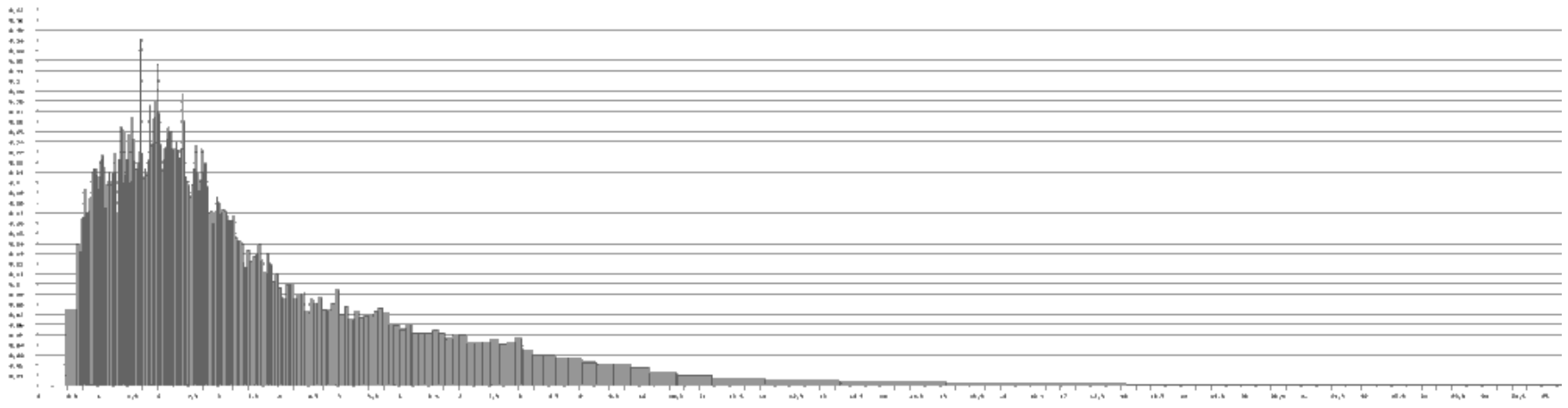

Fig. 10. Bar graph $Z = Y / X$

## 6 Conclusion

While carrying out arithmetic operations over the random variables with unknown distribution, it is advisable to make use of any sample discreteness and build a bar chart of these actions results. At the same time it is essential to ensure that the bar charts in whole as well as each of its half-segments are accurate.

The proposed approach in case with a large number of parameters allows to build a final result bar chart and consider the likelihood of value ranges of a random quantity, rather than to find bulky multi-dimensional convolutions (if it is possible).

Such an approach to assess the results of arithmetic operations allows to select the probable intervals of algebraic equations and their systems we are interested in, differential equations and their systems in case of histogram assessment of empirical probability density distributions for random parameters.

## References

[1] R. E. Moore, Bounding Sets in Function Spaces with Applications to Nonlinear Operator Equations, *SIAM Review*, **20** (1978), 492 – 512. http://dx.doi.org/10.1137/1020068

[2] A. Neumaier, *Interval Methods for Systems of Equations,* Cambridge University Press, 1991. http://dx.doi.org/10.1017/cbo9780511526473

[3] G. Alefeld, J. Herzberger, *Introduction to Interval Computation*, Academic press, 1983.

[4] G. Muscolino, R. Santoro, A. Sofi, Explicit reliability sensitivities of linear structures with interval uncertainties under stationary stochastic excitation, *Structural Safety*, **52** (2015), 219 – 232. http://dx.doi.org/10.1016/j.strusafe.2014.03.001

[5] L. Demidova, Yu. Sokolova, E. Nikulchev, Use of Fuzzy Clustering Algorithms' Ensemble for SVM classifier Development, *International Review on Modelling and Simulations*, **8** (2015), 446 – 457. http://dx.doi.org/10.15866/iremos.v8i4.6825

[6] B.S. Dobronets, O.A. Popova Numerical probabilistic analysis under aleatory and epistemic uncertainty, *Reliable Computing*, **19** (2014), 274 – 289.

[7] P. K. Pollett, The generalized Kolmogorov criterion, *Stochastic Processes and their Applications*, **33** (1989), 29 – 44. http://dx.doi.org/10.1016/0304-4149(89)90064-1

[8] G. Alefeld, G. Mayer, Interval analysis: theory and applications, *Journal of Computational Applied Mathematics*, **121** (2000), 421 – 464. http://dx.doi.org/10.1016/s0377-0427(00)00342-3